\documentclass[10pt,twocolumn]{article}
\usepackage{Nsc2016}
\usepackage[framed]{mcode}

\title{\uppercase{REVISITING HAMMEL ET AL. (1987): DOES THE SHADOWING PROPERTY HOLD FOR MODERN COMPUTERS?}}

\author[1]{\underline{Silva, B. C.}}
\author[2]{Milani, F. L.}
\author[3]{Nepomuceno, E. G.\footnote{Corresponding author. This research has been undertaken in GCOM - Modelling and Control Research Group at Department of Electrical Engineering, Federal University of São João del-Rei, Brazil.}}
\author[4]{Martins, S. A. M.}
\author[5]{Amaral, G. F. V.}

\affil[1]{UFSJ, São João del-Rei, Brazil, bruno.csillva@gmail.com}
\affil[2]{UFSJ, São João del-Rei, Brazil, felipe\_lulli@hotmail.com}
\affil[3]{UFSJ, São João del-Rei, Brazil, nepomuceno@ufsj.edu.br}
\affil[4]{UFSJ, São João del-Rei, Brazil, martins@ufsj.edu.br}
\affil[5]{UFSJ, São João del-Rei, Brazil, amaral@ufsj.edu.br}

\begin{document}
\maketitle

\Abstract{Although computational techniques have been extensively applied in nonlinear science, the reliability evaluation of numerical results does not grow in the same pace. Hammel et al. (Journal of Complexity, 1987, 3(2), 136--145) have proved a theorem that a pseudo-orbit of a logistic map is shadowed by a true orbit within a distance of $10^{-8}$ for $10^{7}$ iterates. We checked this theorem in a modern computational platform, and on the contrary, we found the simulation of the logistic map with less than $10^2$ iterates presents an error greater than $10^{-8}$.}


\section{Introduction}
Numerical computing and Nonlinear Science and Complexity go hand in hand and this has been a long-term relation. In the 60s, it had already been observed the use of computers dedicated to efforts in understanding climate phenomena \citep{Lorenz1963}.  As stated by \citet{Macau2014}, computational techniques are applied to different topics in nonlinear dynamics, such as Synchronization, Bifurcation and Chaos, Complex Networks, Conservative Systems and Nonlinear Partial Differential Equations. \citet{Loz2013} also agrees that computing provides an essential contribution to the analysis of nonlinear dynamical systems.   Thus, beyond using the computer to study systems of such nature, it also becomes important verifying the confiability of its numerical results. Some recent works have explored this issue of checking computer results \citep{Nep2014, Gal2013, NM2016}. Although, it is much more expressive the amount of research simply supported by computers. Added to the theoretical and experimental branches of science, computing is seen as a third one \citep{Ove2001}.

\citet{Hammel1987} were among the first to examine the relation between numerical experiments and the true dynamics of a system \citep{Hayes2006}. By means of the Cray X-MP, a computer costing multi-million dollars at that time, they have reported the shadowing property as valid for a considerable number of iterates when simulating the logistic map for the parameter $r=3.8$ and the initial condition $x_0 =0.4$.  When this case is considered, according to \citet{Hammel1987}, a pseudo-orbit of the logistic map is shadowed by a true orbit within a distance of $10^{-8}$ for $10^{7}$ iterates. Since the theorem is proved for specific conditions, "it can be raised an issue if the computer test constitutes a sufficient condition, and therefore the theorem is proved for all cases, or whether it is a necessary condition, the result may not be valid for all cases", as states \citet{NM2016}.

Many researchers have applied this result in studies concerning the dynamical systems theory since then. More than 100 \citet{Hammel1987} citations have been analyzed through the Scopus online platform and it was observed that, in some cases, the shadowing theorem is seen as a property possible to be generalized for non-hyperbolic systems. \citet{Petrov2013} mention the paper of \citet{Hammel1987} as the one that proves the shadowing property for systems with non-hyperbolic behavior. \citet{Hammel1987} presents the  "first proof of the existence of a shadow for a two-dimensional non-hyperbolic system over a non-trivial length of time", say \citet{Hayes2006}. In addition, \citet {Chaitin-Chatelin1996} also mentions \citet{Hammel1987} when points out considerations about the reliability of numerical results. \citet {Chaitin-Chatelin1996} remarks that the shadowing property had been shown valid for some initial condition and parameter values of the logistic map, without mentioning the cases for which the property fails. 

\section{Purpose}
The aim of this study is to verify whether the shadowing property, as enunciated by \citet{Hammel1987}, is also valid for today's computer systems that adopt the IEEE 754-2008 standard. The same case as the one from that occasion will be studied, namely, the logistic map with $r= 0.4$ and  the initial condition  and $x_0 =3.8$.

\section{Preliminary concepts}
Essential concepts for understanding following sections appear here. The enumerated definitions and the lower bound theorem presented in this section are given by \citet{NM2016}.
\subsection{IEEE 754-2008 standard and the Interval Arithmetic}
The computer works with a limited representation of real numbers, once they are represented by intervals. Therefore, when performing certain mathematical operations, it might often not be possible to store their exact result on the computer memory. The use of round modes, as defined by the IEEE 745-2008 \citep{IEE2008}, is required. Given the combination of applying rounding methods and the way how the computer stores numbers, that is, given the interval arithmetic, mathematical properties such as commutativity, associativity and distributivity, cannot be guaranteed for floating-point operations \citep{Ove2001}. It can be said that the floating-point simulations represent an abstraction of the reality and results from reality may not be easily transferred to the computer \citep{Cor1994}.
\subsection{Orbits and Pseudo-orbits}
Consider the following definitions for orbits and pseudo-orbits.
\begin{defn}
	An orbit is a sequence of values of a map, represented by $x_n = [x_0, x_1, \dots, x_n].$
\end{defn}
\begin{defn}
	Let $i \in \mathbb{N}$ represents a pseudo-orbit, which is defined by an initial condition, an interval extension of $f$, some specific hardware, software and numerical precision standard.
	A pseudo-orbit is an approximation of an orbit and we represent as
	$$ \hat{x}_{i, n} = [\hat{x}_{i, 0}, \hat{x}_{i, 1}, \dots, \hat{x}_{i, n},]$$
	such that, 
	\begin{equation}\label{eq:1}
	|x_n - \hat{x}_{i,n}| \leq \delta_{i,n},
	\end{equation}
	where $\delta_{i,n} \in \mathbb{R}$ is the error and $\delta_{i,n} \geq 0 .$
\end{defn}
\subsection{Shadowing}
Considering the definitions of orbit and pseudo-orbit, \citet{Hammel1987} enunciate a theorem for the logistic map \citep{May1976},
\begin{equation}\label{eq:2}
x_{n+1} = f(x_{n}) = rx_{n}(1-x_n),
\end{equation} 
$n \in \mathbb{N}$ and $a \in \mathbb{R}$.

Through the use of notations for orbit and pseudo-orbit different from those mentioned on the previous definitions, they state that a  pseudo-orbit $ \{p_n\}_{n=0}^{N}$ is shadowed by a true orbit  $ \{x_n\}_{n=0}^{N}$ within shadowing distance $\delta_x = 10^{-8}$, with $p_0 = 0.4$ and $r=3.8$, for $N=10^7$.
\subsection{Lower bound error}
\begin{theo}\label{theo:1}
	Let two pseudo-orbits $\{\hat{x}_{a, n}\}$ and $\{\hat{x}_{b, n}\}$ derived from two interval extensions. Let $\delta_{\alpha, n} = \cfrac{|\hat{x}_{a, n} - \hat{x}_{b, n}|}{2}$ be the lower bound error of a map $f(x)$, then $\delta_{a, n} \geq \delta_{\alpha, n}$ or $\delta_{b, n} \geq \delta_{\alpha, n}$ 
\end{theo}

The proof of this theorem may be found  in \citet{NM2016}.
\section{Methodology}
Consider the following interval extensions of the logistic map (Eq. \ref{eq:2}):
\begin{eqnarray}
    \label{eq:3}
G(X) & = & rX(1-X)\\
\label{eq:4}
H(X) &=& r(X(1-X)).
\end{eqnarray}
It may be observed that both are representations mathematically equivalent \citep{YT2013}.

The computationally generated pseudo-orbits for 100 iterations of \eqref{eq:3} and \eqref{eq:4}, ${\{\hat{x}_{G, n}\}}$ and ${\{\hat{x}_{H, n}\}}$, respectively, are considered for analysis, and, as already mentioned, for 0.4 as initial condition and $r=3.8$. A third pseudo-orbit ${\{\hat{x}_{P, n}\}}$, from a interval extension $P(X)$ of \eqref{eq:2}, was obtained through the use of the variable-precision arithmetic (VPA) Matlab function. 

All the computational tests were performed using the Matlab R2016a on an Intel i5-4440 CPU @ 3.10GHz and Windows 8 (64 bits) operating system.

\section{Results and discussion}
\begin{figure}[htb]
	\centering
	\includegraphics[width=0.9\linewidth]{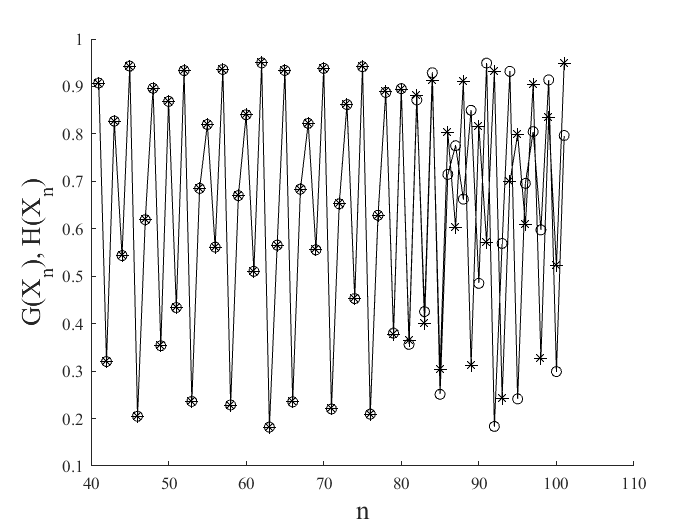}
	\caption{Simulation of $G(X)$ and $H(X)$. {\normalfont Results of $n$ iterates  for $G(X_n)$ (--o--) and $H(X_n)$ (--*--), considering $X_0 = 0.4$ and $r=3.8$.}}
	\label{fig:1}
\end{figure}

\begin{figure}[htb]
	\centering
	\includegraphics[width=0.9\linewidth]{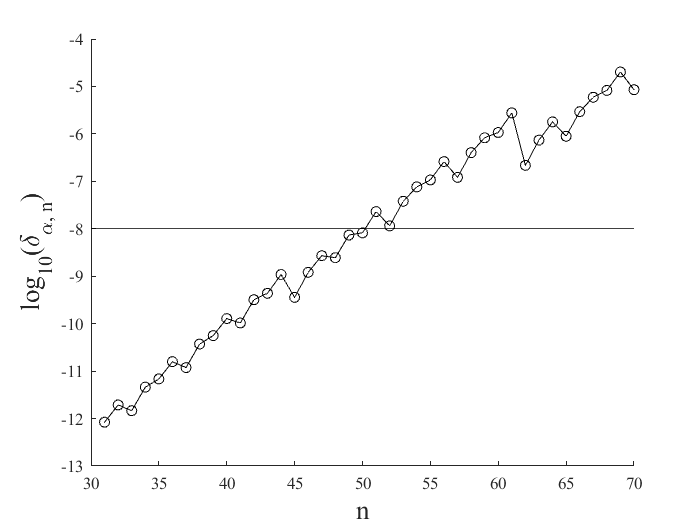}
	\caption{Evolution of the lower bound error $\delta_{\alpha, n}$. {\normalfont Values plotted using the $log_{10}$ for  $\delta_{\alpha, n}$}.}
	\label{fig:2}
\end{figure}

Figure \eqref{fig:1} shows the results for the simulation of interval extensions from Equations \eqref{eq:3} and \eqref{eq:4} after 40 iterates. In the $51^{th}$ iterate it may be already observed a lower bound error $\delta_{\alpha, 51} = 10^{-7.638}$, greater than $10^{-8}$ (Figure \ref{fig:2}). Thus, by Theorem \ref{theo:1}, $\delta_{G, 51} \geq 10^{-7.638}$ or $\delta_{H, 51} \geq 10^{-7.638}$.

\begin{figure}[htb]
	\centering
	\includegraphics[width=0.9\linewidth]{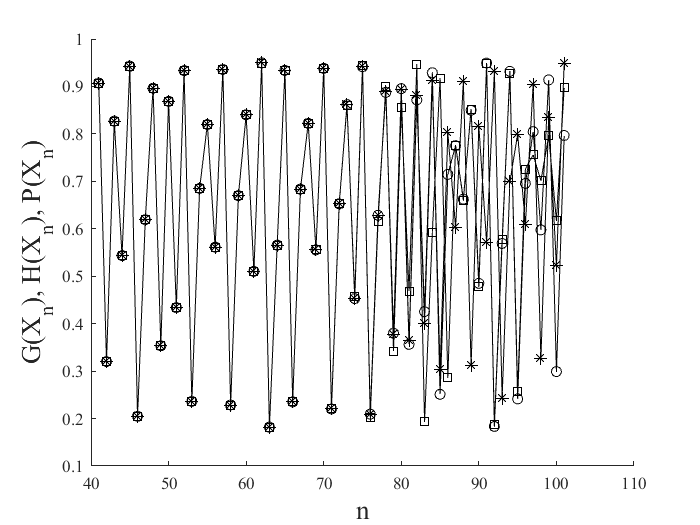}
	\caption{Simulation of $G(X)$, $H(X)$ and $P(X)$. {\normalfont Results of $n$ iterates  for $G(X_n)$ (--o--), $H(X_n)$ (--*--) and $P(X_n)$ (--$\square$--), considering $X_0 = 0.4$ and $r=3.8$.}}
	\label{fig:3}
\end{figure}

In Figure \ref{fig:3}, results for the simulation of $G(X)$ e $H(X)$ are represented, and also, results produced by the iteration process of the logistic map using 1000-digit precision through VPA, $P(X_n)$.

\begin{figure}[htb]
	\centering
	\includegraphics[width=0.9\linewidth]{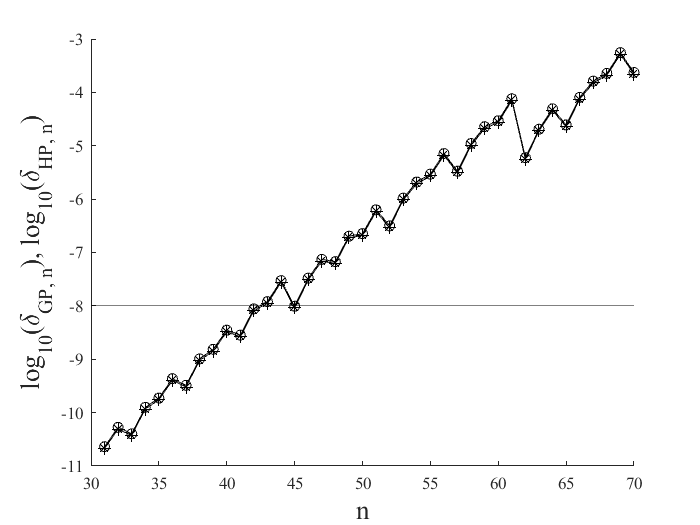}
	\caption{Evolution of $\delta_{GP, n}$ and $\delta_{HP, n}$. {\normalfont Values plotted using the $log_{10}$ for $\delta_{GP, n}$ (--o--) and $\delta_{HP, n}$ (--*--).}}
	\label{fig:4}
\end{figure}

Let $\delta_{GP, n}$ represents $|\hat{x}_{G, n} - \hat{x}_{P, n}|$ and $\delta_{HP, n}$ represents $|\hat{x}_{H, n} - \hat{x}_{P, n}|$. The evolution of $\delta_{GP, n}$ and $\delta_{HP, n}$ is shown in Figure \ref{fig:4}.

As $\{\hat{x}_{P, n}\}$ is produced through symbolic computing of high precision, $\{\hat{x}_{P, n}\}$ is expected to be very close to a true orbit, at least for a number of iterates much lower than the number of iterates, which is the case. Then it is expected that $\delta_{P, n} \approx 0$, $\delta_{G, n} \approx \delta_{GP, n}$, and $\delta_{H, n} \approx \delta_{HP, n}$. Considering that $\delta_{P, 43} \approx 0$, $\delta_{G, 43} \approx 10^{-7.921}$ and $\delta_{H, 43} \approx 10^{-7.954}$, for example.

Therefore, by this logic that combines the Theorem \ref{theo:1} and the use of symbolic computing, it may be said that the results of both interval extensions $G(X)$ and $H(X)$ are far from being good representations of the true orbit of the logistic map for $X_0=0.4$ and $r=3.8$ for $ n>100 $. Thus, the theorem proved by \cite{Hammel1987} does not hold in this computer and software. 

\section{CONCLUSION}
From the presented results and discussion, it may be inferred that the shadowing property, as stated in \cite{Hammel1987} with $r=3.8$ and $x_0=0.4$, cannot be associated to the logistic map simulations using the tested hardware and software that adopts the IEEE-754 floating-point standard.

\section*{ACKNOWLEDGMENTS}
The authors are thankful to the Brazilian agencies CNPq/INERGE, FAPEMIG and CAPES.

\bibliographystyle{unsrtnat}
\bibliography{1603nsc-brunosilva}

\begin{thebibliography}{15}
\providecommand{\natexlab}[1]{#1}
\providecommand{\url}[1]{\texttt{#1}}
\expandafter\ifx\csname urlstyle\endcsname\relax
  \providecommand{\doi}[1]{doi: #1}\else
  \providecommand{\doi}{doi: \begingroup \urlstyle{rm}\Url}\fi

\bibitem[Lorenz(1963)]{Lorenz1963}
Edward~N. Lorenz.
\newblock {Deterministic Nonperiodic Flow}.
\newblock \emph{J. of the Atmospheric Sciences}, 20\penalty0 (2):\penalty0
  130--141, 1963.

\bibitem[Macau and Lambruschini(2014)]{Macau2014}
E.~E.~N. Macau and C.~L.~P. Lambruschini.
\newblock {Advanced computational and experimental techniques in nonlinear
  dynamics}.
\newblock \emph{The European Physical Journal Special Topics}, 223:\penalty0
  2645--2648, 2014.

\bibitem[Lozi(2013)]{Loz2013}
R.~Lozi.
\newblock {Can we trust in numerical computations of chaotic solutions of
  dynamical systems ?}
\newblock In \emph{Topology and dynamics of Chaos, World Scientific Series in
  Nonlinea Science Series A, 84}, pages 63--98, 2013.

\bibitem[Nepomuceno(2014)]{Nep2014}
E.~G. Nepomuceno.
\newblock {Convergence of recursive functions on computers}.
\newblock \emph{The J. of Engineering}, pages 1--3, 2014.

\bibitem[Galias(2013)]{Gal2013}
Z.~Galias.
\newblock {The Dangers of Rounding Errors for Simulations and Analysis of
  Nonlinear Circuits and Systems--and How to Avoid Them}.
\newblock \emph{IEEE Circuits and Systems Magazine}, 13\penalty0 (3):\penalty0
  35--52, 2013.

\bibitem[Nepomuceno and Martins(2016)]{NM2016}
E.~G. Nepomuceno and S.~A.~M. Martins.
\newblock {A lower-bound error for free-run simulation of the polynomial
  NARMAX}.
\newblock \emph{Systems Science {\&} Control Engineering}, pages 1--14, mar
  2016.

\bibitem[Overton(2001)]{Ove2001}
M.~L. Overton.
\newblock \emph{{Numerical Computing with IEEE Floating Point Arithmetic}}.
\newblock SIAM, 2001.

\bibitem[Hammel et~al.(1987)Hammel, Yorke, and Grebogi]{Hammel1987}
S.~M. Hammel, J.~A. Yorke, and C.~Grebogi.
\newblock {Do numerical orbits of chaotic dynamical processes represent true
  orbits?}
\newblock \emph{Journal of Complexity}, 3\penalty0 (2):\penalty0 136--145,
  1987.

\bibitem[Hayes et~al.(2006)Hayes, Jackson, and Young]{Hayes2006}
W.~B. Hayes, K.~R. Jackson, and C.~Young.
\newblock Rigorous high-dimensional shadowing using containment: the general
  case.
\newblock 14\penalty0 (2):\penalty0 329--342, 2006.

\bibitem[Petrov and Pilyugin(2013)]{Petrov2013}
A.~A. Petrov and S.~Y. Pilyugin.
\newblock {Lyapunov functions, shadowing, and topological stability}.
\newblock \emph{arXiv preprint arXiv:1311.3872}, pages 1--11, 2013.

\bibitem[Chaitin-Chatelin(1996)]{Chaitin-Chatelin1996}
F.~Chaitin-Chatelin.
\newblock {Is FINITE precision arithmetic useful for physics?}
\newblock \emph{Journal of Universal Computer Science}, 2\penalty0
  (5):\penalty0 380--395, 1996.

\bibitem[{Institute of Electrical and Electronics Engineers
  (IEEE)}(2008)]{IEE2008}
{Institute of Electrical and Electronics Engineers (IEEE)}.
\newblock \emph{{754-2008 -- IEEE standard for floating-point arithmetic}}.
\newblock IEEE, 2008.

\bibitem[Corless(1994)]{Cor1994}
R.~M. Corless.
\newblock {What good are numerical simulations of chaotic dynamical systems?}
\newblock \emph{Computers {\&} Mathematics with Applications}, 28\penalty0
  (10):\penalty0 107--121, 1994.

\bibitem[May(1976)]{May1976}
R.~M. May.
\newblock {Simple mathematical models with very complicated dynamics}.
\newblock \emph{Nature}, 261\penalty0 (5560):\penalty0 459--67, jun 1976.

\bibitem[Yabuki and Tsuchiya(2013)]{YT2013}
M.~Yabuki and T.~Tsuchiya.
\newblock {Double Precision Computation of the Logistic Map Depends on
  Computational Modes of the Floating-point Processing Unit}.
\newblock \emph{arXiv preprint arXiv:1305.3128}, pages 1--10, 2013.

\end{thebibliography}

\lstset{language=matlab}          

\newpage
\section{Appendix}
\subsection{Matlab code}
\vspace{-0.6cm}	
\begin{lstlisting} 
clear all;
close all;
x=0.4;y=x;
r=3.8;
N=100;
%Pseudo-Orbits
for k=1:N
  x(k+1)=r*x(k)*(1-x(k));
  y(k+1)=r*(y(k)*(1-y(k)));
end
t=1:N+1;
figure(1)
id=41:101;
plot(t(id),x(id),'k-o',t(id),y(id),'-k*');
set(gca, 'FontName', 'Times')
set(gca, 'Box', 'off')
ylabel('G(X_n), H(X_n)', 'fontsize', 16);
xlabel('n', 'fontsize', 16);
%Symbolic - VPA
digits(1000);
z=vpa('4/10');
r=vpa('38/10');
for k=1:N
  z(k+1)=vpa(r.*z(k).*(1-z(k)));
end
figure(2)
plot(t(id),x(id),'k-o',t(id),y(id),'-k*',...
     t(id),z(id),'-ks');
set(gca, 'FontName', 'Times')
set(gca, 'Box', 'off')
ylabel('G(X_n), H(X_n), P(X_n)', 'fontsize', 16);
xlabel('n', 'fontsize', 16);
for(k=1:100)
e(k) = abs(x(k)-y(k))/2;
end
for(k=1:100)
  ex(k) = abs(x(k)-z(k));
end
for(k=1:100)
  ey(k) = abs(y(k)-z(k));
end
figure(3)
ie=31:70;
plot(t(ie),log10(e(ie)),'ko-');
hold on;
plot([30,70],[-8,-8],'k-');
set(gca, 'FontName', 'Times')
set(gca, 'Box', 'off')
ylabel('log_{10}(\delta_{\alpha,n})',...
       'fontsize',16);
xlabel('n', 'fontsize', 16);
figure(4)
plot(t(ie), -8,'-', t(ie),log10(ex(ie)),'ko-',...
     t(ie),log10(ey(ie)),'k*-');
hold on;
plot([30,70],[-8,-8],'k-');
set(gca, 'FontName', 'Times')
set(gca, 'Box', 'off')
s1='log_{10}(\delta_{GP,n}),';
s2='log_{10}(\delta_{HP,n})';
ylabel(strcat(s1,s2),'fontsize', 16);
xlabel('n', 'fontsize', 16);
\end{lstlisting} 

\end{document}